\documentstyle[12pt]{article}
\newcommand{\R}{\mbox{I$\!$R}}

\newcommand{\qed}{{\hfill
{$\rlap{$\sqcap$}\sqcup$}}\\[0.2in]\hspace*{0.5in}}
\newcommand{\qedwh}{{\hfill {$\rlap{$\sqcap$}\sqcup$}}\\[0.2in]}
\newcommand{\bk}{\\[0.03in] \hspace*{0.5in} }

\begin{document}

\vspace*{0.01in}
\begin{center} {\LARGE  Asymptotic Properties of Energy of Harmonic Maps}
\medskip \\ {\LARGE   on Asymptotically Hyperbolic Manifolds}
\medskip \medskip \medskip \medskip\\   {\large Man Chun LEUNG}
\medskip\\
{\large Department of Mathematics, National University of Singapore} \\
{\large
Singapore 119260  \ \ {\tt matlmc@nus.sg}} 
\medskip \medskip \\ 
\end{center}
\begin{abstract}Asymptotic behavior of energy of a harmonic map defined on 
an asymptotically hyperbolic
manifold is considered. Using the growth of energy, we show that a
harmonic
map defined on some asymptotically hyperbolic manifolds has to be constant
if
the total energy is finite, or if the map approaches a point fast enough,
in
terms of a defining function for the boundary.
\end{abstract}
\vspace{0.05in} KEY WORDS \& AMS MS Classifications\footnote{
\hspace*{0.24in}1991 AMS MS Classifications: Primary 58E20; Secondary
58G03}

\vspace{0.4in}

{\bf \Large {\bf \S \,1. \ \ Introduction}}\\

\vspace{0.1in}

Let $(M, g)$ be a complete Riemannian manifold and $(N,
\gamma)$ be a Riemannian manifold. For a smooth map $f : M \to N\,,$ 
the energy of the map $f$ is given by
$$e (f) = g^{ij} \gamma_{\alpha \beta} (f) {{\partial f^{\alpha}}\over
{\partial x_i}} {{\partial f^{\beta}}\over {\partial x_j}}\,.$$
Suppose that $f : M \to N$ is a harmonic
map. A fundamental problem to consider is the relation between the
harmonic map
$f$ and the energy $e (f)$. Burns,
Garber, Ruijsenaars and Seiler [5] show that if $(M, g)$ is the Euclidean
space,
then there are no non-constant harmonic maps with finite energy. As a
special
case of the results obtained in [13], Sealey shows that the same is true
for the hyperbolic space. As often, if the total energy is infinite, then
there
is a growth formula on energy. Price [12] obtains the growth rate on
energy for
harmonic maps defined on the Euclidean space and Karcher and Wood [8]
obtain the
growth rate for simply connected complete Riemannian 
manifolds of bounded negative sectional curvature. The growth
rate is used  by Jin [7] to obtain the following
interesting result.  For $n \ge 3$, if $f : {\R}^n \to (N, \gamma)$ is a
harmonic map such that
$f (x) \to y_o \in N$ as $|x|\to \infty$, then $f$ is a constant map.\bk
In this
paper we study the energy behavior of of a harmonic maps
defined on an asymptotically hyperbolic manifold. Asymptotically
hyperbolic manifolds 
are generalizations of the hyperbolic space with the Poincar\'e metric. 
Other examples of asymptotically hyperbolic manifolds include some
non-compact 
quotients of the hyperbolic space. They are studied by various authors in
connection 
with the Hodge cohomology [9] and Einstein metrics with prescribed
conformal 
structures at the boundary [6]. We consider a harmonic map which is
defined
outside a compact subset of an asymptotically hyperbolic manifold. The
total energy is found to be either infinite or the energy decay
exponentially 
in terms of the distance.  
An example shows that it is possible to have a finite energy harmonic map
defined outside 
a compact subset of a hyperbolic space. \bk
We study harmonic maps defined on the whole asymptotically hyperbolic
manifold. 
It is shown that for a class of complete 
hyperbolic manifolds, there are no non-constant globally
defined harmonic maps with finite energy. The growth rates
of energy for harmonic maps on an asymptotically hyperbolic manifold 
are investigated. We show that if the total energy is infinite, then the
energy grows 
exponentially in terms of the distance. The exponent can be chosen to be
arbitrary close to that of the hyperbolic  space. We note that previous
results 
on growth rate of energy require the manifold to be simply connected.\bk
It is well-known that hyperbolic space supports bounded harmonic maps 
(for example, see [2]). We use the 
asymptotic properties of energy to show that if a harmonic map defined 
on an   
asymptotically hyperbolic manifold approaches a fixed point fast enough,
then the harmonic map has to be 
a constant. More precisely we prove the following result. For an integer
$n \ge
3\,,$  let $M =  {\bf H}^{n + 1}/ \Gamma$ be a
complete hyperbolic manifold without cusps and the exponent of convergence
of the Poincar\'e 
series $\delta (\Gamma ) < n/2 - 1$, or $(M, g)$ is a 
simply connected asymptotically hyperbolic manifold with the sectional 
curvatures of $(M, g)$ satisfying $-b^2 \le K \le - a^2$ with $n a \ge
2b\,.$ 
Let $f : (M, g) \to (N,
\gamma)$ be a $C^2$-harmonic map and $d_N$ be the distance function of
$N$. If
there exist a point $y_o \in N$ and positive constants $C$ and $\sigma$ 
such that  
$$d_N (y_o, f (x)) \le C \rho^{1 + \sigma} (x) \ \ \ \ {\mbox{for}} \ \ 
x \in M\,, \leqno (1.1)$$ 
then $f$ is a constant map. If we assume that $f$ can be extended to a
$C^2$-map on
the closure $\overline M$ (see section 1), then the same conclusion holds
if 
$$d_N (y_o, f (x)) \le C \rho^\sigma (x)\,, \ \ \ \ {\mbox{for}} 
 \ \  x \in M\,.$$
The proof depends on a rather precise 
lower bound on the energy growth of a non-constant harmonic map defined on
an
asymptotically hyperbolic manifold. As in [7], the decay condition (1.1)
gives
an upper bound on the energy growth which contradicts with the
lower bound unless the harmonic map is constant. We apply a regularity 
result to weaken the decay condition.

\vspace{0.4in}

{\bf \Large {\bf \S \,2. \ \ Definitions and Basic Properties}}\\

\vspace{0.1in}

We defined a class of Riemannian manifolds called conformally compact
manifolds. 
Let $M$ be the interior of a compact $(n+1)$-dimensional
manifold $\overline M$ with boundary $\partial M$. Throughout this paper
we 
assume that $n \ge 2$ unless otherwise is stated. A smooth function 
$\rho : \  \overline M \rightarrow R $ is said to be a defining function
for $\partial M$ if 
$$\rho \geq 0\,, \ \ \ \rho^{-1}(0) = \partial M \ \ \ {\mbox{and}} \ \ \ 
d\rho
\not=0 \ \ \ {\mbox{along}} \ \ \ \partial M\,.$$ 
A Riemannian metric $g$
on $M$ is said to be {\it {conformally compact}} if there exist a defining
function $\rho$ for $\partial M$ and a Riemannian metric $h$ on
$\overline M$ such that $g= (1/\rho^2) h$. The Poincar\'e metric
$$\jmath_{ij} = {{4} \over
{({1 - {\mid x \mid}^2})^2}} \  \delta_{ij} \ \ \ \ {\mbox{for}} \ \ 
x \in B^{n + 1}\,, \leqno (2.1)$$
on the unit ball $B^{n + 1}$ is a conformally compact metric, with
defining function 
$$\rho_H (x) = \ 
{{1}\over {2}}(1 - {\mid x \mid}^2) \ \ \ \ {\mbox{for}} \ \ x \in B^{n +
1}\,.$$
\hspace*{0.5in}In general, $(M, g)$ is a complete 
Riemannian manifold which may not be simply connected.  
Mazzeo [9] shows that the sectional curvatures of a conformally 
compact metric approach $-{\mid d\rho
\mid}^2_h$ near the boundary $\partial M$. A conformally compact
metric $g$ is said to be {\it asymptotically hyperbolic} if ${\mid d\rho
\mid}^2_h = 1$ along the boundary. In this case, the sectional
curvatures of an asymptotically hyperbolic metric approach $-1$ near the
boundary. Let $x\in \partial M$ and  $x' = (x_2\,,....,
x_{n+1})$ be a system of local coordinates in a neighborhood of $p$ in
$\partial M$. With $x_1 = \rho\,,$ $(x_1\,, x_2\,,...., x_{n+1}\,)$ is a
coordinate system
in a neighborhood of $p$ in $\overline M$. We may choose a defining
function $\rho$ such that 
in a neighborhood of $\partial M$, the metric $g$ can be written as [6] 
$$g = {{1} \over {\rho^2}} \left[ d\rho^2 +  \sum_{2\le i,\,j \le
n+1} h_{ij}(\rho,\, x')\,dx^idx^j \right]\,. \leqno (2.2)$$
\hspace*{0.5in}In the following we discuss some basic properties we need
about harmonic maps.  
Let $(M, g)$ and $(N, \gamma )$ be Riemannian manifolds of dimension 
$n + 1$ and $m$ respectively. Let $\{y_1,..., y_m\}$ be 
a local coordinate system for $N$. If 
$f : M 
\to N$ is a $C^2$-map, then the energy 
density $e (f)$ of $f$ is defined by 
$$e (f) = g^{ij} \gamma_{\alpha \beta} (f) {{\partial f^{\alpha}}\over
{\partial x_i}} {{\partial f^{\beta}}\over {\partial x_j}}\,. \leqno
(2.3)$$
We observe the  summation convention throughout this paper.  
The map $f$ is called  a harmonic map if $f$ is a critical point of the
energy functional
$$E (f) = \int_{M} e (f) \,\,dg \leqno (2.4)$$ 
with respect to compactly supported variations. 
Let $f$ be a
harmonic map. In local coordinates, $f$ satisfies a system of equations
[3]
$$\Delta f^{\alpha} + g^{ij} \Gamma^{\alpha}_{\mu\nu} {{\partial
f^{\mu}}\over {\partial x_i}} {{\partial f^{\nu}}\over {\partial
x_j}} = 0\,, \ \ \ \alpha = 1, 2,..., m\,, \leqno (2.5)$$ 
where $\Delta$ is the Laplacian operator for $(M, g)$ and
$\Gamma^{\alpha}_{\mu\nu}$ are the  Christoffel symbols of $(N, \gamma)$.
\bk
As consequences of first variation, we have the following two formulas.
The
first one is with respect to a compactly supported vector field $X$ on
$M$. Let
$f : M \to N$ be a harmonic map. We have [12] $$\int_M \{ e (f)
\,{\,{\mbox{div}}\,} X - 2
< df (\bigtriangledown_{e_i} X)\,, \,df (e_i)> \} \,dg = 0\,, \leqno
(2.6)$$
where $\{e_1,..., e_{n +  1}\}$ is an orthonormal basis for $T_p {M'}$ and
$<\ , \ >$ is 
the inner product of $(N, \gamma)$. Assume that $N$ is isometrically
embedded
into ${\bf R}^s$, for some big integer $s$. Let $\vartheta : M \to {\bf
R}^s$
be a compactly supported map, then the first variational formula gives [7]  
$$\int_M g^{\alpha \beta}
\left(2\gamma_{ij} (f) {{\partial f^i}\over {\partial x_\alpha}}
{{\partial
\vartheta^j}\over {\partial x_\beta}} + {{\partial \gamma_{ij} (f)}\over
{\partial y_k}} \vartheta^k {{\partial f^i}\over {\partial x_\alpha}}
{{\partial f^j}\over {\partial x_\beta}} \right) \,dg = 0\,. \leqno
(2.7)$$
The stress-energy tensor $S_f$ of $f$ is defined by 
$$S_f = {1\over 2} e (f) \cdot g - f^* (\gamma)\,. \leqno (2.8)$$
From [3] we have 
$$({\,{\mbox{div}}\,} \,S_f)_i = g^{jk} \bigtriangledown_{e_j} S_{ki} = 0
\ \  
\ \ {\mbox{for}} \ \ i = 1,..., n + 1\,. \leqno (2.9)$$
Let $D$ be a bounded domain in $M$ with 
boundary $\partial D$ and $X$ be a smooth vector field on $D$, then [14]:
$${1\over 2} \int_{\partial D} e (f) <X, {\vec{\bf n}}> ds = \int_D <S_f,
\bigtriangledown X> \,dg + 
\int_{\partial D} <df (X), df ({\vec{\bf n}})> ds\,, \leqno (2.10)$$
where ${\vec{\bf n}}$ is the outward pointing unit normal vector field
along 
$\partial D$ and $ds$ is the induced measure on $\partial D\,.$ Therefore
the 
integral  $\,\int_D < S_f, \bigtriangledown X > \,dg\,$ depends only on
the
behavior of the  vector field along the boundary. This is known as a
conservation law for harmonic map.

\pagebreak

{\bf \Large {\bf \S \,3. \ \ Asymptotic Properties of Energy}}\\

\vspace{0.1in}

In this section we consider harmonic maps defined outside a compact subset
of an 
asymptotically hyperbolic manifold $(M, g)$. We start with an asymptotic
result near the boundary of $M$. Let $K$ be a compact set in $M$, 
which may be empty.  
Given a defining function $\rho$ so that the Riemannian metric $g$ can be
written
in the form as in 2.2. For a positive number $\sigma$, define  
$$D_{\sigma} = \{ x \in M \,\,| \ \rho (x) < \sigma \}\,.$$ 
If $\sigma$ is
small enough, then $D_{\sigma}$ is a diffeomorpic to $\partial M
\times (0, \sigma)\,.$\\[0.2in] 
{\bf Theorem 3.1.} \ \ {\it For an integer $n \ge 3\,,$  
let $f : (M\setminus K\,, g) \to (N\,, h)$ be a
harmonic map. Then either $E (f) = \infty$ or, given any $\epsilon > 0\,,$
there 
exist positive constants $C$ and $\Upsilon$ such that}
$$\int_{D_\sigma} e (f) \,dg \le C \sigma^{\,\ln \,(n - \epsilon)} \ \ \ \
{\mbox{for}} \ \  0 < \sigma < \Upsilon\,.$$
{\it Proof.} \ \ The defining function is chosen in such a way that $g$
can 
be written in the form as in (2.2). Choose a positive number $\sigma_o$
such that 
$ K cap D_{\sigma_o} = \emptyset\,.$ 
For $\sigma > 0$, let 
$$H_\sigma = \{ x \in M \ | \ \rho (x)  = \sigma \}\,.$$ 
There exists a small positive
number $\Upsilon < \sigma_o$ such that for all $0 < \sigma \le \Upsilon$, 
$H_\sigma$ is diffeomorphic to $\partial M$. Assume that $E (f) < \infty$. 
We need to show that
$$\int_\sigma^0 \left(\int_{H_\sigma} e (f) \,dS \right) \,{{d\rho}\over
\rho} = 
\int_\sigma^0 \left(\int_{H_\rho} e (f) \,{{dh}\over
{\rho^n}} \right) \,{{d \rho}\over {\rho}} \ \le \ C \sigma^{\ln ( n -
\epsilon)}
\ \ \ \ {\mbox{for}} \ \ \sigma < \Upsilon\,. \leqno (3.2)$$
Here $\,dS = dh/\rho^n$ is the induced measure on $H_\sigma$ by the
conformally compact metric $g$ and 
$$dh = \sqrt {{\mbox{det}} \,(h_{ij})_{2 \le i, j \le n + 1}} \,\,
dx_2 \cdot \cdot \cdot dx_{n + 1}\,.$$
Denote    
$$D_{\sigma_1\,,\,\sigma_2} =  \{ p \in M \,| \ \sigma_1 <
\rho (p) < \sigma_2 \}\,, \ \ \ 0 < \sigma_1 < \sigma_2 < \Upsilon\,.$$ 
We have $\partial D_{\sigma_1\,,\,\sigma_2} = H_{\sigma_1} \cup
H_{\sigma_2}$.
$E (f) < \infty$ implies that 
$$\int_{\Upsilon}^0 \left(\int_{H_\rho} e (f) \,dS \right) {{d
\rho}\over {\rho}} \ \le \ E (f) < \infty\,.$$
Hence there exists a sequence real numbers $\{
\sigma_i \}_{i = 1}^{\infty}$ such that $\sigma_i < \Upsilon$ and 
$$\lim_{i \to \infty} \sigma_i = 0 \ \ \ \ \ {\mbox{and}} \ \ \ \   
\lim_{i \to \infty} \int_{H_{\sigma_i}} e (f) \,dS =
0\,. \leqno (3.3)$$ 
For any $\sigma < \Upsilon$ and $i$ large, let 
$$X' = - \rho {{\partial}\over {\partial \rho}}$$
be a vector field defined on $D_{\sigma_i\,, \,\sigma}$. Applying (2.10)
we obtain 
\begin{eqnarray*}
& \ & {1\over 2} \int_{\partial D_{\sigma_i\,, \,\sigma}}
e (f) <X', {\vec{\bf n}}> \,dS\\ & = & \int_{D_{\sigma_i\,,\,\sigma}}
<S_f,
\bigtriangledown X'> \,dg + 
\int_{\partial D_{\sigma_i\,,\, \sigma}} <df (X'), df ({\vec{\bf n}})>
\,dS\,.
\end{eqnarray*}
Since 
${\vec{\bf n}} =  \sigma {{\partial}\over {\partial \rho}}$ along
$H_\sigma$ and ${\vec{\bf n}}
= - \sigma_i {{\partial}\over {\partial \rho}}$ along $H_{\sigma_i}$, we
have  
\begin{eqnarray*}
(3.4) \ \ \ \ \ \ \ \ \ \ & \ & {1\over 2}\int_{H_{\sigma_i}} e (f) \,dS  
- {1\over 2}\int_{H_{\sigma}} e (f) \,dS\\
& = &\int_{D_{\sigma_i, \sigma}} 
<S_f, \bigtriangledown X'> \,dg  
 + \int_{H_{\sigma_i}} <df ({\vec{\bf n}}), df ({\vec{\bf n}})> \,dS \ \ \
\
\ \ \ \ \ \ \ \ \ \ \ \ \ \ \ \ \ \ \ \ \ \ \ \ \ \ \ \\\  & \ & \ \ \ \ \
\ -
\int_{H_{\sigma}}  <df ({\vec{\bf n}}), df ({\vec{\bf n}})> \,dS\,.
\end{eqnarray*} As 
$$<df ({\vec{\bf n}}), \,df ({\vec{\bf n}})> \ \ \le \ \ e (f)\,,$$
if we let $i \to \infty$ and apply (3.3), then we have
$$\int_{H_{\sigma}} <df ({\vec{\bf n}}), df ({\vec{\bf n}})> \,dS =
\int_{D_{0, \sigma}} <S_f,
\bigtriangledown X'> \,dg + {1\over 2}\int_{H_{\sigma}} e (f) \,dS\,.
\leqno (3.5)$$ 
Using (2.2) and a calculation show that
$${\,{\mbox{div}}\,} X' =  n - {1\over 2} \,\rho \,h^{kl} {{\partial
h_{kl}}\over {\partial \rho}}\,.$$  Let $x \in M$ with $\rho (x) <
\sigma$. We
may choose local coordinates about $x$ such that 
$$(e_1, e_2,...., e_{n+1}) = (\rho
{{\partial}\over {\partial \rho}}\,, \rho {{\partial}\over {\partial
x_2}}\,,...., 
\rho {{\partial}\over {\partial x_{n+1}}})$$
is an orthonormal basis for $T_x (M)$. Using (2.2), we have
$$\bigtriangledown_{e_1} X' = 0 \ \ \ 
{\mbox{and}} \ \ \ \bigtriangledown_{e_j} X' = e_j - {1\over 2} \, \rho\,
h^{kl} 
\,{{\partial h_{jl}}\over {\partial \rho}} e_k \ \ \ \ {\mbox{for}} \ \ 2
\le j\,, k \le
n + 1\,.$$
From [13], we have 
$$<S_f, \bigtriangledown X'> = {1\over 2} e (f) \,{\,{\mbox{div}}\,} X'
-  <df (e_i), df (e_j)>\, <\bigtriangledown_{e_j} X', e_i>\,.$$
The calculation above shows that  
\begin{eqnarray*}
<S_f, \bigtriangledown X'>
& = & {1\over 2} e (f) \left[ n  - {1\over 2} \rho  h^{kl}  {{\partial
h_{kl}}\over
{\partial \rho}} \right]  + {1\over 2} \rho h^{il} {{\partial
h_{jl}}\over {\partial \rho}} <df (e_i), df (e_j)>\\
& \ & \ \ \ \ \  - \sum_{j = 2}^{n + 1} <df (e_j), df (e_j)>\,.
\end{eqnarray*}
As $\rho \to 0$ on the boundary and $h$ is a smooth metric on the closure 
$\overline M\,,$ all the terms involving $\rho$ can be bounded by
${\delta\over 2} e (f)$ for some 
positive constant $\delta\,.$ Therefore if we choose $\Upsilon$ to be
small,
then we have
$$<S_f, \bigtriangledown X'> 
\ \ge \ {1\over 2} e (f) (n - 2 - \delta) + <df (e_1), df (e_1)> \ \,\ge \
0 
\leqno (3.6)$$
for $\rho (x) < \sigma < \Upsilon\,.$ 
Furthermore $\delta$ can be taken 
to be depending on $\Upsilon$ and can be chosen to be 
as small as possible when $\Upsilon \to 0$. Then (3.5) and (3.6) give
$$\int_{H_{\sigma}} <df ({\vec{\bf n}}), df ({\vec{\bf n}})> \,dS \ \ge \
{1\over 2}\int_{H_{\sigma}} e
(f) \,dS \ \ \ \ {\mbox{for}} \ \  \sigma < \Upsilon\,.\leqno (3.7)$$
Equity holds if and only if 
$$\int_{D_{0, \sigma}} < S_f, \bigtriangledown X' > \,dg = 0\,.$$
From $(3.6)$, this is the case if and only if $e (f) \equiv 0$ on $D_{0,
\,\sigma}$.  By unique continuation, this is equivalent to $f$ being a
constant
map.\bk 
We introduce a change of variables $e^{-t} = \rho$  so that $e^{- \tilde
t}
= \Upsilon$ for some large number $\tilde t\,.$ We have 
$$\int_{\Upsilon}^0 \left(\int_{H_\rho} e (f) \,dS \right)\, {{d
\rho}\over {\rho}} =
\int_{\tilde t}^\infty \left(\int_{H_{e^{-t}}} e (f) \,dS \right) \,dt <
\infty\,.$$
We want to show that for any positive constant $\epsilon\,,$ 
there exists a positive constant $t_o >
\tilde t$  such that for all $t_1 >
t_o$ we have 
$$\int_{t_1}^{t_1 + 1} \left(\int_{H_{e^{-t}}} e (f) \,dS \right) \,dt \
\ge \ (n - 1 -
\epsilon)  \int_{t_1 + 1}^\infty \left(\int_{H_{e^{-t}}} e (f) \,dS
\right) \,dt\,. \leqno
(3.8)$$ 
Assume that (3.8) is not true. Then there exists a positive constant
$\epsilon > 0$ such 
that no matter how large $t_o$ is, there exists a $t_1 > t_o$ such that
$$\int_{t_1}^{t_1 + 1} \left(\int_{H_{e^{-t}}} e (f) \,dS \right) \,dt \ <
\ (n - 1 -
\epsilon)  \int_{t_1 + 1}^\infty \left(\int_{H_{e^{-t}}} e (f) \,dS
\right) \,dt\,. \leqno (3.8')$$
Let $t_2 >
t_1$ and  
$$X =  \chi (t) \,{{\partial}\over {\partial t}}$$
be a vector field defined on $D_{e^{-t_2}, \,e^{-t_1}}$, where $\chi$ is
defined
by 
$$\chi (t) =
\left\{ \begin{array}{ll} 0, & \mbox {\ if $t \le t_1; $}\\
t - t_1, & \mbox{ if $t_1 \ge t \le t_1 + 1\,; $}\\
1,& \mbox {\ if $t_2 \ge t \ge t_1 + 1\,.$}\\
\end{array} \right. $$
A calculation similar to $(3.6)$ gives, for $t \in [t_1, t_1 + 1]$,
\begin{eqnarray*}
<S_f, \bigtriangledown X> & = & {1\over 2} e (f) \left[ n (t - t_1) + 1 -
{1\over 2}
(t -  t_1) \rho \,h^{kl} \,{{\partial h_{kl}}\over {\partial \rho}}
\right]\\
& \ & + {1\over 2}
(t -  t_1) \rho \,h^{kl} \,{{\partial h_{jl}}\over {\partial \rho}} <df
(e_k), df
(e_j)>\\ 
& \ & \ \ \ - \sum_{j = 2}^{n + 1} (t - t_1) <df (e_j), df (e_j)> - <df
(e_1), df
(e_1)>\\
& \ge & {1\over 2} e (f) \left[ (n - 2 - \delta) (t - t_1)\right] + (t -
t_1) <df (e_1), df
(e_1)>\\
& \ &  - <df (e_1), df (e_1)> + {1\over 2} e (f)\\
& \ge & - {1\over 2} e (f)\,.
\end{eqnarray*}
For $t \in [t_1 + 1, t_2]$, we have
\begin{eqnarray*}
<S_f, \bigtriangledown X> & = & {1\over 2} e (f) 
\left[ n - {1\over 2} \rho \,h^{kl} \,{{\partial h_{kl}}\over {\partial
\rho}} 
\right]\\
& \ & \ \ \ + {1\over 2} \rho \,h^{kl} \,{{\partial h_{jl}}\over {\partial
\rho}}
<df (e_k), df (e_j)> - \sum_{j = 2}^{n + 1} <df (e_j), df (e_j)>\\
& \ge & {1\over 2} e (f) [ (n - 2 - \delta) ] + <df (e_1), df (e_1)>\,.
\end{eqnarray*}
Using (2.10) and (3.7) we obtain
\begin{eqnarray*}
(3.9) & \ & {1\over 2} \int_{H_{e^{-t_2}}} e (f) \,dS = \int_{t_1}^{t_1
+ 1} \left(\int_{H_{e^{-t}}} <S_f, \bigtriangledown X> \,dS \right) dt\\
& \ & \ \ \  + \int_{t_1 + 1}^{t_2}
\left(\int_{H_{e^{-t}}} <S_f, \bigtriangledown X> \,dS \right) dt 
+ \int_{H_{e^{-t_2}}} <df (e_1), df (e_1)> \,dS\\ & \ge & - {1\over
2} \int_{t_1}^{t_1 + 1} (\int_{H_{e^{-t}}} e (f) \,dS) dt\\
& \ & \ \ \  + {1\over 2} (n - 1
- \delta) \int_{t_1 + 1}^{t_2} \left(\int_{H_{e^{-t}}} e (f) \,dS \right)
dt 
+ \int_{H_{e^{-t_2}}} <df (e_1), df (e_1)> \,dS\,. \ \ \ \ \ \ \ \ \  \ \
\   \
\ \  
\end{eqnarray*}  
Given any ${\varepsilon}' > 0$, we can choose $t_2$
large enough such that
$$( 1 + {\varepsilon}')  \int_{t_1 + 1}^{t_2} \left(\int_{H_{e^{-t}}} e
(f) \,dS 
\right) dt \ge
\int_{t_1 + 1}^\infty \left(\int_{H_{e^{-t}}} e (f) \,dS \right) dt\,.$$
If we choose $\delta$ and ${\varepsilon}'$ to be small, that is, if we
choose 
$t_o$ to be large, then we have 
$$n - 1  - \epsilon \ \le \ {{n - 1 -
\delta}\over {1 + {\varepsilon}'}}\,.$$ 
($3.8'$) and (3.9) give
\begin{eqnarray*}
& \ & {1\over 2} \int_{H_{e^{-t_2}}} e (f) \,dS\\
& > & -  {1\over 2}
\int_{t_1}^{t_1 + 1} (\int_{H_{e^{-t}}} e (f) \,dS) dt + {1\over 2}
\left({{n - 1 -
\delta}\over {1 + {\varepsilon}'}}\right)  \int_{t_1 + 1}^{\infty}
\left(\int_{H_{e^{-t}}}
e (f) \,dS \right) dt\\ & \ & \ \ \  + \int_{H_{e^{-t_2}}} <df (e_1), df
(e_1)> \,dS\\
& \ge & \int_{H_{e^{-t_2}}} <df (e_1), df (e_2)> \,dS\,.\\ 
\end{eqnarray*} 
The above inequality is in contradiction with (3.7). We note
that by letting $t_o$ to be large, we can choose $\delta$ and
${\varepsilon}'$
to be as small as we like, that is, we can choose $\epsilon$ as small
as we like. Then (3.8) holds for all
$t_1$ large. For
$ t > t_1$ let 
$$F (t) =
\int_t^\infty
\left(\int_{H_{e^{-t}}} e (f) \,dS \right) dt\,.$$ (3.8) implies that 
$$F (t) - F (t + 1) \ge (n - 1 - \epsilon) F ( t + 1)\,,$$
when $t$ is large. That is,
$$F (t) \ge (n - \epsilon) F (t + 1) \ge \cdot \cdot \cdot 
\ge (n - \epsilon)^{m - 1} F (t + m)\,.$$
There exists a positive constant $C$ such that
$$F (t) \le C (n - \epsilon)^{-t} = C e^{- t \ln (n - \epsilon)} 
\ \ \ \ {\mbox{for}} \ t \ {\mbox{large}}\,. \leqno (3.10)$$
Using the relation $e^{-t} = \rho$ we have obtained (3.2). (3.2) or (3.10) 
indicate an exponential decay in terms of the distance to a 
fixed point in $(M, g)$.\qedwh 
{\bf Remark 3.11.} \ \ There are non-constant harmonic maps defined
outside a
compact set of a hyperbolic space which have finite total energy. Let $r =
|x|$, where $x \in B^{n + 1}$. Then 
$$u (r) = \int_r^1 {{(1 - r^2)^{n - 2}}\over {r^{n - 1}}} dt$$
is a non-constant harmonic function defined on $B^{n + 1} \setminus
\{0\}$. 
$f$ is singular at $0$ [1]. The energy density is given by
$$e (u) = \rho_H^2 |{{(1 - t^2)^{n - 2}}\over {t^{n - 1}}}|^2 = O
(\rho_H^{2n - 2})\,.$$
Hence the total energy inside $D_{\rho_H}$ decays 
in order $O (\rho_H^{n - 2})$.\\[0.2in]  
\hspace*{0.5in}Let $M$ be an asymptotically 
hyperbolic manifold. Let $f$ be a harmonic map defined on the {\it whole}
$M$. The
above example demonstrates that we need to
use a global argument in
order to conclude that $f$ is a constant map if $E (f) < \infty$. 
We consider the most interesting case when $M$ is a
hyperbolic manifold. Let $\Gamma$ be a discrete group of isometries on
${\bf
H}^{n + 1}$. Suppose that $\Gamma$ is without parabolic elements and $M =
{\bf
H}^{n + 1}/\Gamma$ is a complete non-compact manifold. It is known that
$M$
is the interior of a compact manifold and the metric is conformally
compact
[11]. Let $\delta (\Gamma )$ be the exponent of convergence of the
Poincar\'e
series corresponding to $\Gamma$. Then for any $\alpha > \delta (\Gamma
)$,
the series
$$\sum_{\gamma \in \Gamma} e^{-\alpha d_H (x, \gamma x)}$$
converges, where $x \in {\bf H}^{n + 1}$ is a fixed point and $d_H$ is the
distance 
function in the hyperbolic space.\\[0.2in]
{\bf Lemma 3.12.} \ \ {\it For an integer $n \ge 3\,,$ let $M = {\bf H}^{n
+
1}/\Gamma$ be a complete hyperbolic manifolds with $\delta (\Gamma ) < n/2
- 1$.
Suppose that
$\Gamma$  has no parabolic elements. Then there is no
non-constant harmonic map on $M$ with finite energy.}\\[0.1in]
{\it Proof.} \ \ Let $\alpha$ be a number such that 
$$\delta (\Gamma ) < \alpha < {n\over 2} - 1\,.$$
We have 
$$\sum_{\gamma \in \Gamma} e^{-\alpha d_H (x, \gamma x)} \ge \sum_{k =
0}^\infty \sharp \{ \gamma \in \Gamma \  | \ k \ < \ d_H (x, \gamma x) \le
k + 1 \}
\,e^{- (k + 1) \alpha}\,. \leqno (3.13)$$ 
The left hand side of (3.13) converges by the comparison test. Therefore
we have
$$\lim_{k \to \infty}  \sharp \{ \gamma \in \Gamma \ | \ k < d_H (x,
\gamma x)
\le k + 1 \} \,e^{- (k + 1) \alpha} = 0\,.$$
Hence there exists a positive constant $C$ such that 
$$\sharp \{ \gamma \in \Gamma \ | \ k < d_H (x, \gamma x) \le k + 1 \} \le
C e^{-
(k + 1) \alpha}\,.$$ 
We obtain
$$\sharp \{ \gamma \in \Gamma | \ 0 < d_H (x, \gamma x) \le k \} \le C
\sum_{i
= 0}^k e^{- (i + 1) \alpha} \le C' \int_0^k e^{\alpha x} dx \le C''
e^{\alpha
k}\,, \leqno (3.14)$$
where $C'$ and $C''$ are positive constants. 
Let $B_x (k) \subset B_x (2k)$ be balls in ${\bf H}^{n + 1}$ with center
at $x$ and 
radius $k$ and $2k$ respectively, where $k = 1, 2,...$. Let $p : {\bf
H}^{n + 1} \to M$ be the covering
map. Assume that there is a harmonic map $f : M \to N$ with finite energy.
Let $C_o$ 
be a positive number such that $E (f) < C_o$. Denote $\tilde f$ the
lifting of $f$
to ${\bf H}^{n + 1}$. For each $y \in B_x (k)$, we have $d_M (p (y), p
(x)) \le
k$, where $d_M$ is the distance in $M$. Hence there exists a fundamental
domain
$\Delta$ such that $y$ and $\gamma x$ are inside $\Delta$ for some $\gamma
\in
\Gamma\,. $ Furthermore $d_H (y, \gamma x) \le k$. By (3.14), the number
of $\gamma x$ inside $B_x (2k)$ is lesser than or equal
to $C'' e^{2\alpha k}$. Therefore the number of fundamental domains needed
to cover $B_x (k)$ is lesser than or equal to $C'' e^{2\alpha k}$. Let
$\tilde E
(\tilde f\,, k)$ be the total energy of $\tilde f$ on $B_x (k)$. Since $E
(f) <
C_o$, we have 
$$\tilde E (\tilde f\,, k) \le C_o C'' e^{2 \alpha k}\,.$$
If $\tilde f$ is a non-constant harmonic map, then $\tilde E (\tilde f\,,
k)$
grows in the order of $e^{(n - 2) k}$ [8]. If $\alpha < n/2 - 1$, then
$\tilde
f$ has to be a constant map. 
Hence $f$ is a constant map.{\hfill {$\rlap{$\sqcap$}\sqcup$}}

\vspace{0.4in}

{\bf \Large {\bf \S \,4. \ \ Energy Growth and Liouville Type Theorem}}\\

\vspace{0.1in}

We show that if $f$ is a harmonic map defined on the whole of $M$ and has
infinite 
energy, then the energy grows exponentially in terms of the
distance.\\[0.2in]
{\bf Theorem 4.1.} \ \ {\it For an integer $n \ge 3\,,$ let $(M, g)$ be an
asymptotically hyperbolic manifold. If $f : (M, g) \to (N, \gamma)$ is a
harmonic map with $E (f) =
\infty$, then for any $\delta > 0$, we have 
$$E (f, \varrho)  = \int_\varrho^\infty \left( \int_{H_\sigma} 
e (f) \,dS \right)\, {{d\rho}\over \rho}   \ \ge \ {C\over 
{\varrho^{n - 2 - \delta} }} \,\,,$$
where $C$ is a positive constant
depending on $f$ and $\delta$.}\\[0.1in]
{\it Proof.} \ \ Let 
$$X = \xi (\rho) \cdot \rho \,{{\partial}\over {\partial \rho}}\,,$$
where $\xi$ is a compactly supported function to be determined later.
A calculation using (2.2) shows that 
$${\,{\mbox{div}}\,} X = \rho {{d\xi}\over {d\rho}} - n\xi +
{1\over 2} \xi \,\rho \,h^{kl} \,{{\partial h_{kl}}\over {\partial
\rho}}\,.$$
Let $x \in M$. We may choose local coordinates above $x$ such that
$$(e_1, e_2,...., e_{n+1}) = (\rho {{\partial}\over {\partial \rho}},
\rho {{\partial}\over {\partial x_2}},...., 
\rho {{\partial}\over {\partial x_{n+1}}})$$
is an orthonormal basis for $T_x (M)$. Then
$$\bigtriangledown_{e_1} X = (\rho \xi') \,e_1 \ \ \ \ {\mbox{and}}
 \ \ \ \ \bigtriangledown_{e_j} X = \left({1\over 2} \xi \,\rho \,h^{kl}
\,{{\partial
h_{jl}}\over {\partial \rho}} \right) \,e_k - \xi e_j$$
for $2 \le j, \,k \le n + 1\,.$ 
Using formula (2.6), we obtain
\pagebreak
\begin{eqnarray*}
(4.2) \ \ \ \ \  0 & = & \int_M e (f) \left\{ \rho \xi' - n\xi + {1\over
2}
\xi \rho h^{kl} {{\partial h_{kl}}\over {\partial \rho}} \right\} \,dg\\
& \ & \  - 2 \int_M \left\{ \rho \xi'
<df (e_1), df (e_1)> - \sum_{j = 2}^{n + 1} \xi <df (e_j), df (e_j)>
\right\} \,dg 
\ \ \ \ \ \ \ \ \ \ \ \ \ \ \ \ \ \ \ \ \ \ \\
& \ & \ - \int_M \rho \,\xi \,h^{kl} \,{{\partial h_{jl}}\over {\partial
\rho}}
<df (e_k), df (e_j)> \,dg\,. 
\end{eqnarray*}
Given any positive number $\delta\,,$ using the fact that 
$$e (f) = \sum_{j = 1}^{n + 1} <df (e_j), df
(e_j)>\,,$$
and (4.2),  we have
\begin{eqnarray*}
(4.3) & \ & - \int_M \rho \,\xi' \,e (f) \,dg + (n - 2 - \delta) \int_M
\xi e
(f) \,dg\\
& = & -\delta \int_M \xi e (f) \,dg  - 2 \int_M \left\{ \rho \xi'
<df (e_1), df (e_1)> \,+ \,\xi <df (e_1), df (e_1)> \right\} \,dg\\
& \ & \ \ \ + \,{1\over 2}
\int_M \xi \,\rho \,h^{kl} \,{{\partial h_{kl}}\over {\partial \rho}} e
(f)
\,dg - \int_M \rho \,\xi \,h^{kl} \,{{\partial h_{jl}}\over {\partial
\rho}} <df (e_k), df
(e_j)> \,dg\,. \ \ \ \ \  \ \ \ \ \ \  \ \ \ 
\end{eqnarray*}
Given $\epsilon > 0\,,$ let  
$$\phi (x) =
\left\{ \begin{array}{ll} 1, & \mbox {\ if $t \le 1; $}\\
1 + {{1-t}\over {\epsilon}}, & \mbox{ if $1 < t < 1 + \epsilon\,; $}\\
0,& \mbox {\ if $t \ge 1 + \epsilon\,.$}\\
\end{array} \right. $$
For $\tau_1 > \tau > 0$, set 
$$\xi (\rho) = \phi \left({{\tau}\over {\rho}}\right) \cdot \left(1 - \phi
\left({{\tau_1}\over {\rho}} \right) \right)\,.$$
Then 
$$\rho {{\partial \xi}\over {\partial \rho}} = - \tau {{\partial
\xi}\over {\partial \tau}} - \tau_1 {{\partial \xi}\over {\partial
\tau_1}}\,. \leqno (4.4)$$
Substitute into (4.3), we obtain
\begin{eqnarray*}
(4.5)&\,&\tau {{\partial}\over {\partial \tau}} 
\left(\int_M \xi e (f) \,dg \right)
+ (n - 2 -\delta) \int_M \xi e (f) \,dg\\
& = & -\delta \int_M \xi e (f) \,dg - \tau_1 {{\partial}\over {\partial
\tau_1}} \left(\int_M \xi e (f) \,dg \right) + 2 \tau {{\partial}\over
{\partial
\tau}} \left(\int_M \xi < df (e_1), df (e_1)> \,dg \right)\\
& \ & \ \ \ \ + 2 \tau_1 {{\partial}\over
{\partial \tau_1}} \left(\int_M \xi < df (e_1), df (e_1)> \,dg \right) 
- 2 \int_M  \xi <df (e_1), df (e_1)> \,dg\\
& \ & \ \ \ \ + {1\over 2}
\int_M \xi \,\rho \,h^{kl} \,{{\partial h_{kl}}\over {\partial \rho}}  e
(f) \,dg - \int_M \rho \,\xi \,h^{kl} \,{{\partial h_{jl}}\over {\partial
\rho}}
<df (e_k), df (e_j)> \,dg\,. \ \ \ \ \ \ \ 
\end{eqnarray*}
Since $h$ is a smooth metric on closure $\overline M$, there exists a
positive
constant $C_o$ such that 
$$|\rho\, h^{kl}\, {{\partial h_{kl}}\over {\partial \rho}}| \le C_o \rho 
\ \ \ \ {\mbox {and}} \ \ \ \ |\rho\, h^{kl} \,{{\partial h_{jl}}\over
{\partial
\rho}} <df (e_k), df (e_j)>| \le C_o \rho e (f)\,.$$
For the positive number $\delta$, if $\tau_1$ and $\tau$ are small, then
the last two 
terms in the right hand side of (4.5) are dominated by the first term. 
Consider the second and third term in the right hand side of (4.5), 
they are negative. As for the fourth term, it is independent of $\tau$. 
If $E (f) = \infty$, then $E (f, \tau) \to \infty$ as 
$\tau \to 0$. Therefore the fourth
term is also dominated by the first term if $\tau$ is small enough.
Finally, the fifth 
term is negative. 
Hence 
$${{\partial}\over {\partial \tau}} \left(\tau^{n - 2 -\delta} \int_M \xi
e (f)
\,dg \right) \le 0\,, \leqno (4.6)$$
which is ture for $\tau$ small. Let $\epsilon \to 0$ and then integrate,
we have
$$E (f, \rho) \ge {{C}\over {\rho^{n - 2 - \delta}}}\,,$$
where $C$ is a positive constant depending on 
the map $f$ and the constant $\delta$.\qed
It is shown in [5] and [12] that if $f : {\bf R}^n \to N$ is a 
non-constant harmonic map, then $E (f)$ cannot be finite. Moreover, [12]
contains 
a growth formula on energy for Euclidean spaces. Later, in [8] and [13],
it is showed that if $f$ is a 
non-constant harmonic map defined on a complete, simply connected manifold
of
suitably pinched negative sectional curvature, then $ E(f) $ again cannot
be finite. Also, 
a growth formula is obtained in [8]. 
Except for the arbitrary positive constant $\delta$, the 
growth formula in [8] agrees with the present formula on hyperbolic
spaces. For further 
discussion on finiteness and growth of energy, we refer to the
comprehensive surveys [3] and [4], 
and the references within.\\[0.2in]  
{\bf Remark 4.7.} \ \ By lemma 3.12 and theorem 4.1, if $M = {\bf H}^{n +
1}/\Gamma$  is a complete hyperbolic manifold without cusps and $\delta
(\Gamma
) < n/2 - 1$, then any non-constant harmonic on $M$ have the energy growth
specified in theorem 4.1.\\[0.2in]  
{\bf Remark 4.8.} \ \ If $(M, g)$ is a {\it simply connected}
asymptotically hyperbolic manifold with negative sectional curvature.
Suppose
that the sectional  curvatures of $(M, g)$ satisfy $-b^2 \le K \le - a^2$
with $n
a
\ge 2b$. Then any harmonic map with finite total energy is a constant map
(cf. [8]).
Hence any non-constant harmonic map on $(M, g)$ have energy growth as in
theorem 4.1.\\[0.2in] 
{\bf Theorem 4.9.} \ \ {\it For an integer $n \ge 3\,,$ let $M$ be a
asymptotically hyperbolic manifold and $f : M \to N$ be a harmonic
map. If there exists a point $y_o \in N$ such that}
$$d_N (f (x), y_o)
\le C
\rho^{1 + \sigma} (x) \ \ \ \ {\mbox{for}} \ \ x \in M\,,$$ 
{\it where $\sigma > 0$ and
$d_N$ is the distant function of  $N$, then $f$ has finite energy. 
In particular, if $M =  {\bf H}^{n + 1}/ \Gamma$ is a
complete hyperbolic manifold without cusps and  $\delta (\Gamma ) < n/2 -
1$, or
$(M, g)$ is a  simply connected asymptotically hyperbolic manifold with
the
sectional  curvatures of $(M, g)$ satisfying $-b^2 \le K \le - a^2$ with
$n a \ge
2b$, then $f$ is a constant map.}\\[0.1in]
{\it Proof.} \ \ Let $\{y_1,..., y_m\}$ be a local coordinate system above
a neighborhood $U$ 
of $y_o$ such that $y_o = 0$ and
$${{\partial \gamma_{ij} (y)}\over {\partial y_k}} y_k + 2 \gamma_{ij}
(y) \ge \gamma_{ij} (y) \ \ \ \ {\mbox{for}} \ \ y\in U\,. \leqno (4.10)$$
See [7] for a discussion on the existence of such a local coordinate
neighborhood.     
In such a local coordinate system we have
$$|f^i (x)| = |y_i \circ f (x)| \le C
\rho^{1 + \sigma} (x) \ \ \ \ {\mbox{for}} \ \ x \in M\,.$$
Near $\partial M$, let $\{x_1 = \rho, \,x_2,...., x_{n + 1} \}$ be 
local coordinates for $M\,,$ where $\rho$ is a defining function such that 
the Riemannian metric $g$ can be 
expressed in the form as in (2.2). Then
$$E (f) = \int_M \rho^2 h^{\alpha \beta} \gamma_{ij} {{\partial
f^i}\over {\partial x_\alpha}} {{\partial f^j}\over {\partial x_\beta}}
{{1}\over {\rho^{n + 1}}} \,dh\,.$$
Let $\xi$ be defined as in the proof of theorem 4.1. The first variational
formula (2.7) with respect to
$$\vartheta (x) = \xi (\rho (x)) f$$ 
gives (cf. [7])
\begin{eqnarray*}
(4.11) \ \ \ \ \ \ \ \ \ \ & \ & \int_M \left[{{\partial \gamma_{ij}
(f)}\over
{\partial y_k}} f^k + 2 \gamma_{ij} (f) \right] h^{\alpha \beta}
{{\partial f^i}\over
{\partial x_\alpha}} {{\partial f^j}\over {\partial x_\beta}} \xi (\rho
(x))
{1\over {\rho^{n -1} }} \,dh\\ & = & -2 \int_M h^{\alpha \beta}
\gamma_{ij}
{{\partial f^i}\over {\partial x_\alpha}} f^j {{\partial \xi}\over
{\partial
x_\beta}}  {1\over {\rho^{n -1}}} \,dh\,. \ \ \ \ \ \ \ \ \ \ \ \ \ \ \ \
\ \ \ \
\  \ \ \ \ \ \ \ \ \ \ \ \ \ \ \ \ \ \ \ \ \ \ \ \ \ \ \ \ \ \ 
\end{eqnarray*}
We have
$${{\partial \xi}\over {\partial x_\beta}} = {1\over {\epsilon \tau}}
\cdot {{\tau^2}\over {\rho^2}} {{\partial \rho}\over {\partial
x_\beta}} \ \ \ \ {\mbox{for}} \ \ (1 + \epsilon) \tau > \rho >
\tau\,. \leqno (4.12)$$
If we let $\epsilon \to 0$, we obtain
\begin{eqnarray*}
(4.13) \ \ \ \ \ \ & \ & \int_{M (\tau) \setminus M (\tau_1)}
\left[{{\partial
\gamma_{ij} (f)}\over {\partial y_k}} f^k + 2 \gamma_{ij} (f) \right]
\,h^{\alpha
\beta} {{\partial f^i}\over {\partial x_\alpha}} {{\partial f^j}\over
{\partial
x_\beta}} {1\over {\rho^{n -1} }} \,dh + R (\tau_1)\\  
& = & -2 \int_{\partial M (\tau)}  h^{\alpha \beta} \gamma_{ij} {{\partial
f^i}\over {\partial x_\alpha}} f^j {{\partial \rho}\over {\partial
x_\beta}}  {1\over {\rho^{n -1} }} \,ds\,, \ \ \ \ \ \ \ \ \ \ \ \ \ \ \ \
\ \ \ 
\ \ \ \ \ \ \ \ \ \ \ \ \ \ \ \ \ \ \ \ \ \ \ \ \ \ \ 
\end{eqnarray*}
where $\tau_1 > \tau$ and $R (\tau_1)$ is a term depending on $\tau_1$ but
not on
$\tau$. Here 
$ds$ is the measure of $\partial M_\tau$ induced by the Riemannian metric 
$h$. Recall that 
$$M_\tau = \{\, x \in M \ | \ \rho (x) > \tau\,\}\,.$$
Let
$$F (\tau) = \int_{M (\tau) \setminus M (\tau_1)}  \gamma_{ij} (f)
h^{\alpha \beta} {{\partial f^i}\over {\partial x_\alpha}} {{\partial
f^j}\over {\partial x_\beta}} {1\over {\rho^{n -1} }} \,dh + R (\tau_1)\,. 
\leqno (4.14)$$ 
Differentiating both side of (4.14) we  have
$${{d F}\over {d\tau}} = - \int_{\partial M (\tau)} \gamma_{ij} (f))
h^{\alpha \beta} {{\partial f^i}\over {\partial x_\alpha}} {{\partial
f^j}\over {\partial x_\beta}} {1\over {\rho^{n -1}}} \,ds\,.$$
Using H\"older's inequality and the fact that $\rho$ and
$h^{\alpha \beta}$ are smooth functions on the closure $\overline M$, we
obtain 
\begin{eqnarray*}
(4.15) \ \ \ \ \ \ \ & \ & \left(- \int_{\partial M (\tau)}  h^{\alpha
\beta}
\gamma_{ij} {{\partial f^i}\over {\partial x_\alpha}} f^j {{\partial
\rho}\over
{\partial x_\beta}}  {{1}\over {\rho^{n -1}}} \,ds \right)^2 \\
& \le & C' \, \left(\int_{\partial M (\tau)} \gamma_{ij} (f))
h^{\alpha \beta} {{\partial f^i}\over {\partial x_\alpha}} {{\partial
f^j}\over {\partial x_\beta}} {1\over {\rho^{n -1} }} \,ds \right) \times
\left(\int_{\partial M (\tau)} \gamma_{ij} f^i f^j {1\over {\rho^{n -1} }}
\,ds\right)\,, 
\ \ \ \ \ \ \ \ 
\end{eqnarray*}
where $C'$ is a positive constant depending on $\rho$ and $h$. 
Since $d_N (f (x), y_o) = O (\rho^{1 + \sigma} (x))$ for some $\sigma >
0$, 
we have
$$\int_{\partial M (\tau)} h_{ij} f^i f^j {1\over {\rho^{n -1}}} \,ds
\le {{C''}\over {\tau^{n - 3 - 2\sigma}}}\,, \leqno (4.16)$$
for some constant $C'' > 0$. Using (4.14), (4.15) and (4.16), we obtain
$$F^2 (\tau) \le - {{C''}\over {\tau^{n - 3 - 2\sigma}}} {{d F}\over
{d\tau}}\,,$$
or
$${{d}\over {d\tau}} \left({1\over {F (\tau)}} \right) \ge {{\tau^{n - 3 -
\sigma}}\over {C''}}\,. \leqno (4.17)$$ 
Assume that the harmonic map $f$ has infinite energy.  
Then $F (\tau) \to \infty$ as $\tau \to 0$. Hence
$${1\over {F (\tau)}} = \int^\tau_0 {{d}\over {d\tau}} \left({1\over {F
(\tau)}} \right) \,d\tau \ge {1\over {C''}} \int^\tau_0 \tau^{n - 3 -
2\sigma}
\,d\tau\,. $$ Or 
$$E (f, \tau ) \le C'' \tau^{-(n - 2 - 2\sigma)}\,. \leqno (4.18)$$
As the total energy is infinity, we can apply theorem 4.1. If we choose
the
positive constant 
$\delta$ in theorem 4.1 to be lesser than $2 \sigma$, then we have
obtained a contradiction. Hence the total energy of $f$ is finite. In case 
$M =  {\bf H}^{n + 1}/ \Gamma$ is a
complete hyperbolic manifold without cusps and  $\delta (\Gamma ) < n/2
-1$, 
or $(M, g)$ is a 
simply connected asymptotically hyperbolic manifold with the sectional 
curvatures of $(M, g)$ satisfying $-b^2 \le K \le - a^2$ with $n a \ge
2b$,
then $f$ is a constant map by remark 4.7 and 4.8. \qed
We introduce normal coordinates $\{ r, \omega_2\,,....
\omega_m \}$ about $y_o$, where $r$ is the radial distance to $q$ and
$(\omega_2\,,....
\omega_m) \in S^{m-1}$. Denote 
$$f^1 = r \circ f, \ f^2 = \omega_2 \circ f,..., \ f^m = \omega_m \circ
f\,.$$ 

{\bf Lemma 4.19.} \ \ {\it Let $M$ be an asymptotically
hyperbolic manifold and $f : M \to
N$ be a harmonic map. Assume that 
$f^i, 2 \le i \le m$, can be extended as $C^2$-functions on 
$\overline M\,.$ 
If there exists a point $y_o \in N$ such that} 
$$d_N (f (x), y_o) \le C \rho^{\nu} (p) \ \ \ {\rm  for \ \ some} \ \  \nu
>
0\\,,.$$ 
{\it then} 
$$d_N (f (p)\,, q) = O (\rho^n)\,.$$
{\it Proof.} \ \ As above,  $\{x_1 = \rho, x_2,...., x_{n + 1} \}$ is a
local coordinate system for $M$. The equations for harmonic maps are of
the form
$$J (f^1) = \Delta f^1 + \rho^2 h^{\alpha \beta} \Gamma^{1}_{ij}
{{\partial f^{i}}\over {\partial x_\alpha}} {{\partial f^{j}}\over
{\partial x_\beta}} = 0\,. \leqno (4.20)$$
We regard this as an equation on $f^1$. The terms involving $f^i$
and its derivatives, for $2 \le i \le m$, are considered to be
coefficients of the equation. By the assumption we have $f^1 = O
(\rho^\nu)$. Let $L$ be the linearization of $J$ about the zero
function, that is,
$$L (v) = \Delta v + \rho^2 h^{\alpha \beta} \left(\Gamma^{1}_{1j} 
{{\partial v}\over {\partial x_\alpha}} {{\partial f^{j} }\over
{\partial x_\beta}} + \Gamma^{1}_{i1} 
{{\partial f^i}\over {\partial x_\alpha}} {{\partial v }\over
{\partial x_\beta}} \right)\,, \ \ \ 2 \le i,\, j \le m\,. \leqno (4.21)$$
The normal operator [10] of $L$ is $\Delta$, the Laplacian operator for
$(M, g)$, 
which has $0$ and $n$ as the indicial roots. The argument in [10]
shows that $f^1 = O (\rho^n)$.\qed
We can apply lemma 4.19 to weaken the decay assumption in theorem 4.9. 

\pagebreak

\end{document}